\documentclass [10pt] {article}

 \usepackage [cp1251] {inputenc}
 \usepackage [english] {babel}
 \usepackage {indentfirst}
 \usepackage {hyperref}
 \usepackage [curve] {xypic}
 \usepackage {tikz}
 \usepackage {mathtools}
 \usepackage {amsthm}
 \usepackage {amssymb}
 \usepackage [mathscr] {eucal}

 \DeclareFontEncoding {LS1}{}{}
 \DeclareFontSubstitution {LS1}{stix2}{m}{n}
 \DeclareSymbolFont {stix2-symbols}{LS1}{stix2scr}{m}{n}

 \let \til\~     

 \def \?{{?}}

 \def \O{{\displaystyle\raisebox{0.99ex}{$\scriptscriptstyle\boldsymbol<$}\mkern-4.1mu{|}}}
 \def \o{{\scriptstyle\raisebox{0.51ex}{$\scriptscriptstyle<$}\mkern-4.4mu{|}}}
 \def \0{{\mathchoice\O\O\o\o}}

 \DeclareMathSymbol \sbighash{\mathord}{stix2-symbols}{"9F}      
 \def \shash{{\mathchoice
 {\vcenter{\hbox{$\scriptstyle\sbighash$}}}
 {\vcenter{\hbox{$\scriptstyle\sbighash$}}}
 {\vcenter{\hbox{$\scriptscriptstyle\sbighash$}}}
 {\vcenter{\hbox{$\scriptscriptstyle\sbighash$}}}
 }}
 \def \uc{\shash}                

 \def \AA{{\mathscr A}}          

 \def \<{\langle}                
 \def \>{\rangle}

 \def \LQ{$\scriptscriptstyle<$}                 
 \def \`{{\mathchoice
 {\raisebox{0.20ex}\LQ}
 {\raisebox{0.20ex}\LQ}
 {\raisebox{0.10ex}\LQ}
 {\raisebox{0.05ex}\LQ}
 }}

 \def \RQ{$\scriptscriptstyle>$}                 
 \def \'{{\mathchoice
 {\raisebox{0.20ex}\RQ}
 {\raisebox{0.20ex}\RQ}
 {\raisebox{0.10ex}\RQ}
 {\raisebox{0.05ex}\RQ}
 }}

 \newcommand* {\bxd} [5] {                       
 \vcenter {
 \vspace{#4}
 \hbox{\hspace{#3}\tikz{
 \node [inner sep=0ex,minimum size=#1,line width=#2,draw] {$#5$}
 }\hspace{#3}}
 \vspace{#4}
 }
 }

 \def \sbigeCop{{\mathchoice                     
 {\bxd{3.1ex}{0.15ex}{0.15ex}{0.23ex}{\displaystyle\boldsymbol\sqcup}}
 {\bxd{2.2ex}{0.12ex}{0.13ex}{0ex}{\scriptstyle\boldsymbol\sqcup}}
 {\bxd{1.5ex}{0.10ex}{0.13ex}{0ex}{\scriptscriptstyle\boldsymbol\sqcup}}
 {\bxd{1.1ex}{0.07ex}{0.09ex}{0ex}{\scriptscriptstyle\sqcup}}
 }}
 \def \bigeCop{\mathop\sbigeCop}

 \DeclareMathOperator \ord{ord}
 \DeclareMathOperator \supp{supp}
 \DeclareMathOperator \dist{dist}
 \DeclareMathOperator \diam{diam}

 \def \hsm{\leftthreetimes}      

 \def \shs{{\mathchoice         
 {{\scriptstyle\leftthreetimes}}
 {{\scriptstyle\leftthreetimes}}
 {{\scriptscriptstyle\leftthreetimes}}
 {{\scriptscriptstyle\leftthreetimes}}
 }}
 \def \hs{\mathbin\shs}

 \def \le{\leqslant}
 \def \ge{\geqslant}

 \def \xto{\xrightarrow}

 \def \-{\overline}              
 \def \~{\widetilde}
 \def \^{\widehat}

 \renewcommand* {\%} [2] {       
 \overset{#2}#1
 }

 \renewcommand* {\$} [3]
 {\overset{#2}{\underset{#3}#1}}

 \newcommand* {\head} [1] {
 \subsubsection* {#1}
 }

 \newenvironment* {claim} [1] []
 {\begin{trivlist}\item [\hskip\labelsep {\bf #1}] \it}
 {\end{trivlist} }

 \newenvironment* {demo} [1] []
 {\begin{trivlist}\item [\hskip\labelsep {\it #1}] }
 {\end{trivlist} }

\begin {document}

 \title {\large\bf
         Homotopy similarity of maps.
         Compositions}

 \author {\normalsize\rm
          S.~S.~Podkorytov}

 \date {}

 \maketitle

 \vspace {-2\bigskipamount}
 
 \begin {abstract} \noindent
 We describe the behaviour
 of
 the homotopy similarity relations
 and
 finite-order invariants
 under the function $[X,Y]\to[X,Z]$
 induced by a map $Y\to Z$
 strongly $r$-similar to the constant map.
 \end {abstract}


 \head {\S~1. Introduction}


 This paper continues \cite{sim} and \cite{sim-2-str}.
 We adopt notation and conventions thereof.
 Let
 $X$,
 $Y$,
 and
 $Z$
 be cellular spaces,
 $X$
 and
 $Y$
 compact.
 We prove the following two theorems.

 \begin {claim} [1.1. Theorem.]
 Let
 maps $a,a':X\to Y$ satisfy $a\%\sim{p-1}a'$
 and
 a map $b:Y\to Z$ satisfy $b\%\approx{q-1}\0^Y_Z$
 ($p,q\ge1$).
 Then
 the maps $b\circ a,b\circ a':X\to Z$ satisfy
 $$
 b\circ a\%\sim{pq-1}b\circ a'.
 $$
 \end {claim}

 \begin {claim} [1.2. Theorem.]
 Let
 $b:Y\to Z$ be a map such that
 $b\%\approx{q-1}\0^Y_Z$
 ($q\ge1$).
 Let
 $L$ be an abelian group
 and
 $h:[X,Z]\to L$ be a homotopy invariant.
 Then
 the invariant
 $$
 f:[X,Y]\to L,
 \qquad
 [a]\mapsto h([b\circ a]),
 $$
 satisfies
 $$
 q\ord f\le\ord h.
 $$
 \end {claim}

 We conjecture that
 the assumption of strong $(q-1)$-similarity
 in these statements
 can be replaced by that of $(q-1)$-similarity.


 \head {\S~2. Coherence of an ensemble of compositions}


 Let
 $E$ be a nonempty finite set
 and
 continuous functions
 $\phi_e:X\to[0,1]$,
 $e\in E$,
 form a partition of unity:
 $$
 \sum_{e\in E}
 \phi_e
 =
 1.
 $$
 Consider the unbased map
 $$
 \phi=(\phi_e)_{e\in E}:X\to\Delta E.
 $$
 Introduce
 the function
 $$
 \theta^\phi:
 Y^X\times(Z^Y)^{(\Delta E)}
 \to
 Z^X,
 \qquad
 \theta^\phi(a,s):
 X
 \xto{x\mapsto\phi(x)\hs a(x)}
 \Delta E\hsm Y
 \xto{\uc^Y(s)}
 Z
 $$
 [
 equivalently,
 $\theta^\phi(a,s):x\mapsto s(\phi(x))(a(x))$
 ],
 and
 the homomorphism
 $$
 (\theta^\phi):
 \<Y^X\>\otimes\<(Z^Y)^{(\Delta E)}\>
 \to
 \<Z^X\>,
 \qquad
 \`a\'\otimes\`s\'
 \mapsto
 \`\theta^\phi(a,s)\'.
 $$

 We fix a metric on $X$.
 A number $\lambda>0$ is called
 a {\it Lebesgue number\/} of an open cover $\Gamma$ of $X$ if
 any set $V\subseteq X$ of diameter at most $\lambda$
 is contained in some $G\in\Gamma$
 (according to \cite[proof of theorem~3.3.14]{S}).

 \begin {claim} [2.1. Proposition.]
 Let $\Gamma$ be an open cover of $X$ with
 a Lebesgue number $\lambda$.
 Suppose that
 $$
 \diam\supp\phi_e\le\epsilon,
 \qquad
 e\in E,
 $$
 where
 $\epsilon>0$
 and
 $(pq-2)\epsilon\le\lambda$
 ($p,q\ge1$).
 Let an ensemble $A\in\<Y^X\>$ satisfy
 $$
 A\$={p-1}\Gamma0.
 $$
 Let $S\in\<(Z^Y)^{(\Delta E)}\>$ be a fissile ensemble such that
 \begin {equation} \label {fil}
 \`\Xi^{\Delta E}(\0^Y_Z)\'-S
 \in
 \<(Z^Y)^{(\Delta E)}\>^{(q)}_Y.
 \end {equation}
 Then
 $$
 (\theta^\phi)(A\otimes S)
 \in
 \<Z^X\>^{(pq)}.
 $$
 \end {claim}

 \begin {demo} [Proof.]
 Take a set $V\subseteq X$ with $|V|\le pq-1$.
 We should show that
 $$(\theta^\phi)(A\otimes S)|_V=0
 $$
 in $\<Z^{(V)}\>$.
 Consider the equivalence on $V$ generated by
 all pairs $(x_1,x_2)$ with $\dist(x_1,x_2)\le\epsilon$.
 Let $V_i$,
 $i\in I$,
 be the classes of this equivalence.
 Clearly,
 $\diam V_i\le(pq-2)\epsilon$
 and
 $\dist(V_i,V_j)>\epsilon$
 for $i\ne j$.
 For $i\in I$,
 put
 $$
 F_i
 =
 \{\,e\mid\phi_e|_{V_i}\ne0\,\}
 \subseteq
 E.
 $$
 Clearly,
 $\phi(V_i)\subseteq\Delta F_i$.
 We have $F_i\ne\varnothing$ because
 $V_i\ne\varnothing$.
 We have $F_i\cap F_j=\varnothing$
 for $i\ne j$
 because
 $\diam\supp\phi_e\le\epsilon$
 and
 $\dist(V_i,V_j)>\epsilon$.
 Thus
 we have the layout
 $$
 F_*
 =
 \{\,F_i\mid i\in I\,\}
 \in
 \AA(E).
 $$

 Consider the homomorphisms
 $$
 \rho_1:
 \<Y^X\>
 \to
 \bigotimes_{i\in I}
 \<Y^{(V_i)}\>,
 \qquad
 \`a\'
 \mapsto
 \bigotimes_{i\in I}
 \`a|_{V_i}\',
 $$
 and
 $$
 \rho_2:
 \<(Z^Y)^{(\Delta E)}\>
 \to
 \bigotimes_{i\in I}
 \<(Z^Y)^{(\Delta F_i)}\>,
 \qquad
 \`s\'
 \mapsto
 \bigotimes_{i\in I}
 \`s|_{\Delta F_i}\'.
 $$
 From now on,
 let decorated $\rho$'s denote similar homomorphisms.
 For $i\in I$,
 we have,
 similarly to $\theta$ and $\theta^\phi$,
 the function
 $$
 \theta^\phi_i:
 Y^{(V_i)}\times(Z^Y)^{(\Delta F_i)}
 \to
 Z^{(V_i)},
 \qquad
 \theta^\phi(d,t):
 x
 \mapsto
 t(\phi(x))(d(x)),
 $$
 and
 the homomorphism
 $$
 (\theta^\phi_i):
 \<Y^{(V_i)}\>\otimes\<(Z^Y)^{(\Delta F_i)}\>
 \to
 \<Z^{(V_i)}\>,
 \qquad
 \`d\'\otimes\`t\'
 \mapsto
 \`\theta^\phi_i(d,t)\'.
 $$
 Clearly,
 \begin {equation} \label {res}
 \theta^\phi(a,s)|_{V_i}
 =
 \theta^\phi_i(a|_{V_i},s|_{\Delta F_i}),
 \qquad
 a\in Y^X,
 \
 s\in(Z^Y)^{(\Delta E)}.
 \end {equation}

 We have
 \begin {equation} \label {rho2S}
 \rho_2(S)
 =
 \bigotimes_{i\in I}
 S|_{\Delta F_i}.
 \end {equation}
 Indeed,
 consider the commutative diagram with sendings
 $$
 \xymatrix {
 \<(Z^Y)^{(\Delta E)}\>
 \ar[r]^-{\?|_{\Delta[F_*]}}
 \ar[d]_-{\rho_2}
 &
 \<(Z^Y)^{(\Delta[F_*])}\>
 \ar[dl]^-{\-\rho_2}
 \\
 \bigotimes\limits_{i\in I}
 \<(Z^Y)^{(\Delta F_i)}\>,
 &
 }
 \qquad
 \xymatrix {
 {\scriptstyle
 S
 }
 \ar@{|->}[r]
 \ar@{|->}[d]
 &
 \bigeCop\limits_{i\in I}
 {\scriptstyle
 S|_{\Delta F_i}
 }
 \ar@{|->}[dl]
 \\
 {\scriptstyle
 \bigotimes\limits_{i\in I}
 S|_{\Delta F_i}
 }
 &
 }
 $$
 The horizontal sending holds because
 $S$ is fissile.
 The diagonal one is obvious.
 The vertical sending,
 which is \eqref{rho2S},
 follows.

 We have the commutative diagram with sendings
 $$
 \xymatrix @R1.4ex {
 &
 {
 \scriptstyle
 \rho_1(A)\otimes\rho_2(S)
 }
 &
 \\
 &
 \bigotimes\limits_{i\in I}
 \<Y^{(V_i)}\>
 \hfill
 \otimes
 \bigotimes\limits_{i\in I}
 \<(Z^Y)^{(\Delta F_i)}\>
 \ar@{=}[dd]
 &
 \\
 {
 \scriptstyle
 A\otimes S
 }
 \ar@{|->}@/^3ex/[uur]
 & &
 {
 \scriptstyle
 \rho_1(A)
 }
 \ar@{|->}@/_3ex/[uul]
 \\
 \<Y^X\>\otimes\<(Z^Y)^{(\Delta E)}\>
 \ar[uur]^-{\rho_1\otimes\rho_2}
 \ar[dd]_-{(\theta^\phi)}
 &
 \bigotimes\limits_{i\in I}
 (
 \<Y^{(V_i)}\>
 \otimes
 \<(Z^Y)^{(\Delta F_i)}\>
 )
 \ar[dd]_-{
 \bigotimes\limits_{i\in I}
 (\theta^\phi_i)
 }
 &
 \bigotimes\limits_{i\in I}
 \<Y^{(V_i)}\>
 \ar[uul]_-{
 \?
 \otimes
 \bigotimes\limits_{i\in I}
 S|_{\Delta F_i}
 }
 \ar[ddl]^*!/r4ex/{\labelstyle
 \bigotimes\limits_{i\in I}
 (\theta^\phi_i)(\?\otimes S|_{\Delta F_i})
 =
 \bigotimes\limits_{i\in I}
 h_i
 }
 \\ \\
 \<Z^X\>
 \ar[r]^-{\rho_3}
 \ar[dd]_-{\?|_V}
 &
 \bigotimes\limits_{i\in I}
 \<Z^{(V_i)}\>,
 &
 \\ \\
 \<Z^{(V)}\>
 \ar[uur]^-{\-\rho_3}_-{\cong}
 & &
 }
 $$
 where
 $$
 h_i=(\theta^\phi_i)(\?\otimes S|_{\Delta F_i}):
 \<Y^{(V_i)}\>
 \xto{\?\otimes S|_{\Delta F_i}}
 \<Y^{(V_i)}\>
 \otimes
 \<(Z^Y)^{(\Delta F_i)}\>
 \xto{(\theta^\phi_i)}
 \<Z^{(V_i)}\>.
 $$
 Commutativity of the upper-left pentagon follows from \eqref{res}.
 Commutativity of the rest is obvious.
 Clearly,
 $\-\rho_3$ is an isomorphism.
 The first sending is obvious.
 The second one follows from \eqref{rho2S}.
 We should show that $A\otimes S$ goes to zero
 under the composition in the left column.
 By the diagram,
 it suffices to show that
 \begin {equation} \label {rho1A}
 \bigl(
 \bigotimes_{i\in I}
 h_i
 \bigr)
 (\rho_1(A))
 =0.
 \end {equation}

 Let $J\subseteq I$ consist of those $i$
 for which $|V_i|\ge q$.
 We have $|J|\le p-1$ because
 $|V|\le pq-1$.

 Take
 $i\in I\setminus J$
 and
 $d\in Y^{(V_i)}$.
 We show that
 \begin {equation} \label {hid}
 h_i(\`d\')=\`\0^{V_i}_Z\',
 \end {equation}
 where $\0^{V_i}_Z\in Z^{(V_i)}$ is
 (of course)
 the unbased map
 that takes the whole $V_i$ to the basepoint $\0_Z$.
 Consider the unbased maps
 $$
 D:V_i\to\Delta E\hsm Y,
 \qquad
 x
 \mapsto
 \phi(x)\hs d(x),
 $$
 and
 $D'=D|_{V_i\to D(V_i)}$.
 We have the commutative diagram with sendings
 $$
 \xymatrix @C1.4ex{
 {\scriptstyle
 S
 }
 \ar@{|->}[d]
 &
 {\scriptstyle
 \`\Xi^{\Delta E}(\0^Y_Z)\'
 }
 \ar@{|->}[d]
 &
 \<(Z^Y)^{(\Delta E)}\>
 \ar[rrrr]^-{\<\uc^Y\>}
 \ar[d]^-{(\theta^\phi_i)(\`d\'\otimes\?|_{\Delta F_i})}
 &&&&
 \<Z^{(\Delta E\hsm Y)}\>
 \ar[d]^-{\?|_{D(V_i)}}
 \\
 {\scriptstyle
 h_i(\`d\')
 }
 &
 {\scriptstyle
 \`\0^{V_i}_Z\'
 }
 &
 \<Z^{(V_i)}\>
 &&&&
 \<Z^{(D(V_i))}\>.
 \ar[llll]_-{\<Z^{(D')}\>}
 }
 $$
 Commutativity is checked directly.
 The first sending holds by the definition of $h_i$.
 The second one is obvious.
 Consider the difference $R=\`\Xi^{\Delta E}(\0^Y_Z)\'-S$.
 By \eqref{fil},
 $\<\uc^Y\>(R)\in\<(Z^Y)^{\Delta E\hsm Y}\>^{(q)}$.
 Since $|D(V_i)|\le|V_i|\le q-1$,
 we have $\<\uc^Y\>(R)|_{D(V_i)}=0$.
 The equality \eqref{hid} follows
 by the diagram.

 We have the commutative diagram
 $$
 \xymatrix {
 \<Y^X\>
 \ar[rr]^-{\rho_1}
 \ar[drr]_-{\rho_1'}
 &&
 \bigotimes\limits_{i\in I}
 \<Y^{(V_i)}\>
 \ar[rr]^-{
 \bigotimes\limits_{i\in I}
 h_i
 }
 \ar[d]^-{\pi}
 &&
 \bigotimes\limits_{i\in I}
 \<Z^{(V_i)}\>
 \\
 &&
 \bigotimes\limits_{i\in J}
 \<Y^{(V_i)}\>
 \ar[rr]^-{
 \bigotimes\limits_{i\in J}
 h_i
 }
 &&
 \bigotimes\limits_{i\in J}
 \<Z^{(V_i)}\>,
 \ar[u]_-{\sigma}
 }
 $$
 where $\pi$ and $\sigma$ are defined by the rules
 \begin {alignat*} {3}
 \pi:
 &
 \bigotimes_{i\in I}
 \`d_i\'
 \mapsto
 \bigotimes_{i\in J}
 \`d_i\',
 &\qquad
 &
 \textrm{
 $d_i\in Y^{(V_i)}$
 ($i\in I$),
 }
 \\
 \sigma:
 &
 \bigotimes_{i\in J}
 C_i
 \mapsto
 \bigotimes_{i\in I}
 C_i,
 &\qquad
 &
 \textrm{
 $C_i\in\<Z^{(V_i)}\>$
 ($i\in I$),
 $C_i=\`\0^{V_i}_Z\'$
 for $i\notin J$.
 }
 \end {alignat*}
 Commutativity of the square follows from \eqref{hid}.
 We show that
 $\rho_1'(A)=0$.
 By the diagram,
 \eqref{rho1A} will follow.

 For each $i\in I$,
 there is $G_i\in\Gamma$ such that
 $V_i\subseteq G_i$
 because
 $\diam V_i\le(pq-2)\epsilon\le\lambda$.
 Put
 $$
 H=
 \{\0_X\}
 \cup
 \bigcup_{i\in J}
 G_i
 \subseteq
 X.
 $$
 Since $|J|\le p-1$,
 $H\in\Gamma(p-1)$.
 We have the commutative diagram
 $$
 \xymatrix {
 \<Y^X\>
 \ar[rr]^-{\rho'_1}
 \ar[d]_-{\?|_H}
 &&
 \bigotimes\limits_{i\in J}
 \<Y^{(V_i)}\>
 \\
 \<Y^H\>.
 \ar[urr]_-{\^\rho_1'}
 &&
 }
 $$
 Since $A\$={p-1}\Gamma0$,
 we have $A|_{H}=0$.
 By the diagram,
 $\rho_1'(A)=0$.
 \qed
 \end {demo}


 \head {\S~3. Exploiting Proposition~2.1}


 \begin {claim} [3.1. Corollary.]
 Let
 an ensemble $A\in\<Y^X\>$,
 $$
 A=
 \sum_i
 u_i\`a_i\',
 $$
 satisfy $A\%={p-1}0$
 and
 a map $b:Y\to Z$ satisfy $b\%\approx{q-1}\0^Y_Z$
 ($p,q\ge1$).
 Then
 there exists an ensemble $C\in\<Z^X\>$,
 \begin {equation} \label {C}
 C=
 \sum_{i,j}
 u_iv_j\`c_{ij}\',
 \end {equation}
 where
 $c_{ij}\sim b\circ a_i$
 and
 \begin {equation} \label {aff}
 \sum_j
 v_j
 =1,
 \end {equation}
 such that
 $C\%={pq-1}0$.
 \end {claim}

 \begin {demo} [Proof.]
 By \cite[Corollary~6.2]{sim},
 there is an ensemble $\~A\in\<Y^X\>$,
 $$
 \~A=
 \sum_i
 u_i\`\~a_i\',
 $$
 where $\~a_i\sim a_i$,
 such that
 $\~A\$={p-1}\Gamma0$
 for some open cover $\Gamma$ of $X$.
 Let $\lambda$ be a Lebesgue number of $\Gamma$.
 Choose
 $\epsilon>0$ such that
 $(pq-2)\epsilon\le\lambda$
 and
 a partition of unity
 $$
 \sum_{e\in E}
 \phi_e
 =
 1,
 $$
 where
 $E$ is a nonempty finite set
 and
 $\phi_e:X\to[0,1]$ are continuous functions such that
 $$
 \diam\supp\phi_e\le\epsilon,
 \qquad
 e\in E.
 $$
 Form the unbased map
 $$
 \phi=(\phi_e)_{e\in E}:X\to\Delta E.
 $$
 Since $b\%\approx{q-1}\0^Y_Z$,
 there is a fissile ensemble
 $S\in\<(Z^Y_b)^{(\Delta E)}\>\subseteq\<(Z^Y)^{(\Delta E)}\>$,
 $$
 S=
 \sum_j
 v_j\`s_j\',
 $$
 such that
 $\Xi^{\Delta E}(\0^Y_Z)-S\in\<(Z^Y)^{(\Delta E)}\>^{(q)}_Y$.
 Put
 $$
 C=(\theta^\phi)(\~A\otimes S)
 $$
 (see \S~2 for $(\theta^\phi)$).
 By Proposition~2.1,
 $C\%={pq-1}0$.
 We have
 $$
 C=
 \sum_{i,j}
 u_iv_j\`\theta^\phi(\~a_i,s_j)\'.
 $$
 Since
 $\~a_i\sim a_i$
 and
 $s_j\in(Z^Y_b)^{(\Delta E)}$,
 we have $\theta^\phi(\~a_i,s_j)\sim b\circ a_i$.
 The equality \eqref{aff} holds because
 $S$ is
 fissile
 and
 thus
 affine.
 \qed
 \end {demo}

 \begin {demo} [Proof of Theorem~1.1.]
 Since $a\%\sim{p-1}a'$,
 we have
 $$
 \sum_i
 u_i\`a_i\'
 \%={p-1}
 \`a'\'
 $$
 in $\`Y^X\'$,
 where $a_i\sim a$.
 By Corollary~3.1,
 $$
 \sum_{i,j}
 u_iv_j\`c_{ij}\'
 \%={pq-1}
 \sum_j
 v_j\`c'_j\'
 $$
 in $\<Z^Y\>$,
 where
 $c_{ij}\sim b\circ a_i\sim b\circ a$,
 $c'_j\sim b\circ a'$,
 and
 the equality \eqref{aff} holds.
 By \cite[Theorem~7.3]{sim},
 $b\circ a\%\sim{pq-1}b\circ a'$.
 \qed
 \end {demo}

 \begin {demo} [Proof of Theorem~1.2.]
 Suppose that
 $\ord h\le pq-1$
 for some integer $p\ge 1$.
 It suffices to show that
 $\ord f\le p-1$.
 Take an ensemble $A\in\<Y^X\>$,
 $$
 A=
 \sum_i
 u_i\`a_i\',
 $$
 such that
 $A\%={p-1}0$.
 Corollary~3.1 yields an ensemble $C\in\<Z^X\>$
 of the form \eqref{C}
 with
 $c_{ij}\sim b\circ a_i$
 and
 $v_j$ satisfying \eqref{aff}
 such that
 $C\%={pq-1}0$.
 We have
 $$
 \sum_i
 u_if([a_i])
 =
 \sum_{i,j}
 u_iv_jh([b\circ a_i])
 =
 \sum_{i,j}
 u_iv_jh([c_{ij}])
 \overset{(*)}=
 0,
 $$
 where $(*)$ holds because
 $\ord h\le pq-1$
 and
 $C\%={pq-1}0$.
 Thus
 $\ord f\le p-1$.
 \qed
 \end {demo}



 \noindent
 \href{mailto:ssp@pdmi.ras.ru}{\tt ssp@pdmi.ras.ru}

 \noindent
 \url{http://www.pdmi.ras.ru/~ssp}

 \end {document}